\documentclass[12pt]{article}
\usepackage{amsmath}
\usepackage{amsfonts}
\usepackage{amssymb}
\usepackage{longtable}

\newtheorem{theorem}{Theorem}
\newtheorem{lemma}{Lemma}

\def\tr{\mathop\text{\rm tr}\nolimits}
\def\Id{\mathop\text{\rm id}\nolimits}
\def\id{\mathop\text{\rm id}\nolimits}
\def\Ric{\mathop{{\rm Ric}}\nolimits}
\def\tRic{\mathop{\widetilde{\rm Ric}}\nolimits}

\begin{document}

\title{Conformally flat Lorentzian manifolds with special holonomy groups}
\author{Anton S. Galaev}

\maketitle

\begin{abstract} The local classification of  conformally flat Lorentzian
manifolds with special holonomy groups is obtained. The
corresponding local metrics  are certain extensions of Riemannian
spaces of constant sectional curvature to Walker metrics.

\end{abstract}

\section{Introduction and the main results}

 It is known \cite{Kur} that
a conformally flat Riemannian manifold is either a product of two
spaces of constant sectional curvature, or it is a product of a
space of constant sectional curvature with an interval, or its
restricted holonomy group is the identity component of the
orthogonal group. The last condition represents the general case
and among various manifolds satisfying the last condition one can
emphasize only the spaces of constant sectional curvature.

One says that the connected holonomy group of an
indecomposable pseudo-Riemannian manifold is special if it is  different from the connected component of the
pseudo-orthogonal group \cite{Bryant2}. In this case the holonomy
group defines a special geometry on the manifold. For example, a
pseudo-Riemannian manifold of signature $(r,s)$ is
pseudo-K\"ahlerian if and only if its holonomy group is contained
in ${\rm U}(\frac{r}{2},\frac{s}{2})$. We see that there are no
conformally flat Riemannian manifolds with special holonomy
groups.

In the case of pseudo-Riemannian manifolds, the holonomy group can
be weakly irreducible, this means that it does not preserve any
non-degenerate proper vector subspace of the tangent space, and
not irreducible in the same time, i.e. it may preserve a
degenerate vector subspace of the tangent space.

 The main result of the present paper is the complete local
description of conformally flat Lorentzian manifolds $(M,g)$ with
weakly irreducible not irreducible holonomy groups, which are the
only special holonomy groups of Lorentzian manifolds. Let $\dim
M=n+2\geq 4$. The holonomy algebra, i.e. the Lie algebra of the
holonomy group, $\mathfrak{g}\subset\mathfrak{so}(1,n+1)$ of such
manifold preserves an isotropic line of the tangent space, which
is identified with the Minkowski space $\mathbb{R}^{1,n+1}$. Hence
$\mathfrak{g}$ is contained in the maximal subalgebra of
$\mathfrak{so}(1,n+1)$ preserving an isotropic line. This algebra
is denoted by $\mathfrak{sim}(n)$ and it admits the decomposition
$$\mathfrak{sim}(n)=(\mathbb{R}\oplus\mathfrak{so}(n))\ltimes\mathbb{R}^n.$$ The classification of
the Lorentzian holonomy algebras
$\mathfrak{g}\subset\mathfrak{sim}(n)$ can be found in \cite{ESI}.
Any manifold $(M,g)$ with such holonomy algebra (locally) admits a
distribution $\ell$ of isotropic lines. Such manifolds are called
the Walker manifolds \cite{WalkerBook}. On any such manifold
$(M,g)$ there exist local coordinates $v,x^1,...,x^n,u$ such that
the metric $g$ has the form
\begin{equation}\label{Walker} g=2dvdu+h+2Adu+H (d
u)^2,\end{equation} where $h=h_{ij}(x^1,...,x^n,u)d x^id x^j$ is
an $u$-dependent family of Riemannian metrics,  $A=A_i(x^1, \ldots
, x^n,u)d x^i$ is an $u$-dependent family of one-forms, and $H$ is
a local function on $M$. The vector field $\partial_v$ defines the
parallel distribution of isotropic lines. An important class of
Walker manifolds form pp-waves that are given locally by
\eqref{Walker} with $A=0$, $h=\sum_{i=1}^n(dx^i)^2$, and
$\partial_vH=0$, see \cite{L-N}. The pp-waves are exactly Walker
manifolds with the commutative holonomy algebra
$\mathfrak{g}\subset\mathbb{R}^n\subset\mathfrak{sim}(n)$.

We use the Einstein summation convention for the indeces
$i,j,k=1,...,n$; we use the denotation $\dot f=\partial_u f$ for
any function $f$.

On a Walker manifold $(M,g)$  we define the canonical function
$\lambda$ from the equality
$$\Ric p=\lambda p,$$ where $p$ is any local vector field tangent to the distribution $\ell$, and $\Ric$ is the Ricci operator.
If the metric $g$ is written in the form \eqref{Walker}, then $\lambda=\frac{1}{2}\partial^2_vH$, and the scalar curvature of $g$ satisfies $$s=2\lambda+s_0,$$ where $s_0$ is the scalar curvature of $h$. The form of a conformally flat Walker matric will depend on the vanishing of the function $\lambda$. In the general case we obtain the following result.

\begin{theorem}\label{ThC2} Let $(M,g)$ be a conformally flat (i.e. with zero
Weyl curvature tensor) Walker Lorentzian manifold of dimension
$n+2\geq 4$. Then in a neighborhood of each point of $M$ there
exist coordinates $v,x^1,...,x^n,u$ such that
$$g=2dvdu+\Psi\sum_{i=1}^n(dx^i)^2+2Adu+(\lambda(u)v^2+vH_1+H_0)(d
u)^2,$$ where
\begin{align*}\Psi&=\frac{4}{\left(1-\lambda(u)\sum_{k=1}^n(x^k)^2\right)^2},\\
A&=A_idx^i,\quad
A_i=\Psi\left(-4C_k(u)x^kx^i+2C_i(u)\sum_{k=1}^n(x^k)^2\right),\\
H_1&=-4C_k(u)x^k\sqrt{\Psi}-\partial_u\ln\Psi+K(u),\\
H_0(x^1,...,x^n,u)&=\left\{\begin{array}{ll}\frac{4}{\lambda^2(u)}\Psi\sum_{k=1}^nC^2_k(u)+\sqrt{\Psi}\left(a(u)\sum_{k=1}^n(x^k)^2+D_k(u)x^k+D(u)\right),&\textrm{ if } \lambda(u)\neq 0,\\
16\left(\sum_{k=1}^n(x^k)^2\right)^2\sum_{k=1}^nC^2_k(u)+\tilde a(u)\sum_{k=1}^n(x^k)^2+\tilde D_k(u)x^k+\tilde D(u),& \textrm{ if } \lambda(u)=0\end{array}\right.
\end{align*} for some functions $\lambda(u)$, $a(u)$, $\tilde a(u)$, $C_i(u)$, $D_i(u)$,
$D(u)$, $\tilde D_i(u)$, $\tilde D(u)$.

The scalar curvature of $g$ equals to $-(n-2)(n+1)\lambda(u)$.
\end{theorem}

If the function $\lambda$ is locally constantly zero, or it is non-vanishing, then the above metric may be simplified.

 \begin{theorem}\label{ThC2A} Let $(M,g)$ be a conformally flat Walker Lorentzian manifold\\ of dimension $n+2\geq 4$.
\begin{itemize}
\item[1)] If the function $\lambda$ is non-vanishing at a point, then
in a neighborhood of this point there exist coordinates
$v,x^1,...,x^n,u$ such that
$$g=2dvdu+\Psi\sum_{i=1}^n(dx^i)^2+(\lambda(u)v^2+vH_1+H_0)(d
u)^2,$$ where
\begin{align*}\Psi&=\frac{4}{\left(1-\lambda(u)\sum_{k=1}^n(x^k)^2\right)^2},\\
H_1&=-\partial_u\ln\Psi,\quad
H_0=\sqrt{\Psi}\left(a(u)\sum_{k=1}^n(x^k)^2+D_k(u)x^k+D(u)\right).
\end{align*}
\item[2)] If $\lambda\equiv 0$ in a neighborhood of a point, then
in a neighborhood of this point there exist coordinates
$v,x^1,...,x^n,u$ such that
$$g=2dvdu+\sum_{i=1}^n(dx^i)^2+2Adu+(vH_1+H_0)(d
u)^2,$$ where
\begin{align*}
A&=A_idx^i,\quad
A_i=C_i(u)\sum_{k=1}^n(x^k)^2,\quad H_1=-2C_k(u)x^k\\
H_0&=\sum_{k=1}^n(x^k)^2\left(\frac{1}{4}\sum_{k=1}^n(x^k)^2\sum_{k=1}^nC^2_k(u)-(C_k(u)x^k)^2+\dot C_k(u)x^k+a(u)\right)+D_k(u)x^k+D(u).
\end{align*}
In particular, if all $C_i\equiv 0$, then the metric can be rewritten in the form
\begin{equation}\label{cfppw}
 g=2dvdu+\sum_{i=1}^n(dx^i)^2+a(u)\sum_{k=1}^n(x^k)^2(d
u)^2.
\end{equation}
\end{itemize}
\end{theorem}

Thus Theorem \ref{ThC2A} gives the local form of a conformally
flat  Walker metric in the neighborhoods of any point where
$\lambda$ is non-zero or constantly zero. Such points represent a
dense subset of the manifold. Theorem \ref{ThC2} describes also
the metric in the neighborhoods of points at that the function
$\lambda$ vanishes, but it is not locally zero, i.g. in the
neighborhoods of isolated zero points of $\lambda$.

Next, we find the holonomy algebras of the obtained metrics and
check when the metrics are decomposable.

\begin{theorem}\label{ThC3}
Let $(M,g)$ be as in Theorem \ref{ThC2}.
\begin{itemize}
\item[{\bf 1)}] The manifold $(M,g)$ is locally indecomposable if and only if
there exists a coordinate system with one of the properties:
\begin{itemize} \item[$\bullet$] $\dot\lambda\not\equiv 0$,
 \item[$\bullet$] $\dot\lambda\equiv 0$, $\lambda\neq 0$, i.e. $g$ can be written as in part 1) of Theorem \ref{ThC2A},\\ and $\sum_{k=1}^nD^2_k+(a+\lambda D)^2\not\equiv 0$,
 \item[$\bullet$] $\lambda\equiv 0$, i.e. $g$  can be written as in part 2) of Theorem \ref{ThC2A},
 and $\sum_{k=1}^nC^2_k+a^2\neq 0$.
 \end{itemize}
Otherwise, the metric can be written in the form
$$g=\Psi\sum_{k=1}^n(dx^k)^2+2dvdu+\lambda v^2(du)^2,\quad \lambda\in\mathbb{R}.$$
The holonomy algebra of this metric is trivial if and only if
$\lambda=0$. Otherwise, the holonomy algebra is isomorphic to
$\mathfrak{so}(n)\oplus\mathfrak{so}(1,1)$.

\item[{\bf 2)}] Suppose that $(M,g)$ is  locally indecomposable. Then its holonomy algebra is isomorphic to $\mathbb{R}^n\subset\mathfrak{sim}(n)$ if and only if $\lambda^2(u)+\sum_{k=1}^nC^2_k(u)\equiv 0$  for all coordinate systems.
In this case $(M,g)$ is a pp-wave, and $g$ is given by
\eqref{cfppw}. Otherwise, the holonomy algebra is isomorphic to
$\mathfrak{sim}(n)$.
\end{itemize}\end{theorem}

In Section \ref{secCdec}, the result of Kurita \cite{Kur} is
extended to the case of pseudo-Riemannian manifolds. In Section
\ref{secCWeyl}, an expression for the curvature tensor and the
Weyl conformal curvature tensor $W$ for a Walker metric is given.
It seems that these expressions are given here for the first time.
Sections \ref{secCdokvo}, \ref{secpT2A} and \ref{secpT3} are
dedicated to the proofs of the main theorems. We rewrite the
equation $W=0$ as a system of partial differential equations and
find appropriate systems of coordinates such that the complete
solution of this system can be found. In Section \ref{secCRic},
the Ricci operator for the obtained metrics  is computed.

In Section \ref{secCdim4}, the case of dimension 4 is considered.
Possible holonomy algebras of conformally flat 4-dimensional
Lorentzian manifolds are classified  also in \cite{H-L}. The
metric \eqref{cfppw} in dimension 4 is given in \cite{St}. In
\cite{H-L}, it is posed the problem to construct an example of
 conformally flat metric with the holonomy algebra
$\mathfrak{sim}(2)$ (which is denoted in \cite{H-L} by $R_{14}$).
An attempt to construct such metric was done in \cite{GT}. We show
that the metric constructed there is in fact decomposable and its
holonomy algebra is $\mathfrak{so}(1,1)\oplus\mathfrak{so}(2)$.
 Thus in this paper we get conformally flat metrics with the
holonomy algebra $\mathfrak{sim}(n)$ for the first time, and even
more,  we find all such metrics.

The field equations of Nordstr\"om's theory of gravitation, which
appeared before Einstein's theory, are the following:
$$W=0,\quad s=0,$$
see \cite{Nor,Rav,Vizgin}. All metrics from Theorem \ref{ThC2} in
dimension 4 and metrics from part 2) of Theorem \ref{ThC2A} in
bigger dimensions
 provide examples of solutions of these equations. Thus we have found all solutions to
Nordstr\"om's gravity with holonomy algebras contained in
$\mathfrak{sim}(n)$. Similarly, the Einstein equation on
Lorentzian manifolds with such holonomy algebras was studied in
\cite{G-P}. In this case it is impossible to obtain the complete
solution, but the examples of solutions have interesting physical
interpretations \cite{G-P}.

An important fact is that a simply connected
conformally flat spin Lorentzian manifold  admits the space of
conformal Killing spinors of maximal dimension \cite{Baum08}.

It would be interesting to obtain examples of conformally flat
Lorentzian manifolds satisfying some global geometric properties,
e.g. important are globally hyperbolic Lorentzian manifolds with
special holonomy groups \cite{Baz}.

The projective equivalence of 4-dimensional conformally flat
Lorentzian metrics with special holonomy algebras was studied
recently in \cite{Hall12}. There are many interesting works about
conformally flat (pseudo-)Riemannian, and in particular Lorentzian
manifolds. Let us mention the works
\cite{AleKim78,Kir92,Sla94,E-I,Hon,H-K,H-K07}.

The results of this paper are used in \cite{2sym} for the
classification of Lorentzian manifolds satisfying the condition
$\nabla^2R=0$.

\section{Decomposability of conformally flat pseudo-Riemannian
manifolds}\label{secCdec}

 In \cite{Kur}, Kurita proved
the following theorem for the case of  Riemannian manifolds.

\begin{theorem} Let $(M,g)$ be an $n$-dimensional conformally flat
Riemannian manifold. Then its local restricted holonomy group
$H_x$ ($x\in M$) is in general ${\rm SO}(n)$. If $H_x\neq {\rm
SO}(n)$, then for some coordinate neighborhood $U$ of $x$ one of
the following holds: \begin{itemize} \item[1)] $H_x$ is identity
and the metric is flat in $U$;
\item[2)] $H_x={\rm SO}(k)\times {\rm SO}(n-k)$ and $U$ is a direct product
of a $k$-dimensional manifold of constant sectional curvature $K$
and an $(n-k)$-dimensional manifold of constant sectional
curvature $-K$ ($K\neq 0$);
\item[3)] $H_x={\rm SO}(n-1)$ and $U$ is a direct product of a straight
line (or a segment) and an $(n-1)$-dimensional manifold of
constant sectional curvature.\end{itemize} \end{theorem}

We generalize this theorem for the case of pseudo-Riemannian
manifolds. We also make it more precise.

\begin{theorem}\label{ThC1} Let $(M,g)$ be a conformally flat
pseudo-Riemannian manifold of signature $(r,s)$ with the
restricted holonomy group ${\rm Hol}^0(M,g)$. If $(M,g)$ is not
flat, then one of the following holds: \begin{itemize} \item[1)]
${\rm Hol}^0(M,g)={\rm SO}(r,s)$;
\item[2)]
${\rm Hol}^0(M,g)$ is weakly irreducible and not irreducible (in
particular, it preserves a degenerate subspace of the tangent
space);
\item[3)] ${\rm Hol}^0(M,g)={\rm SO}(r_1,s_1)\times {\rm SO}(r-r_1,s-s_1)$
and each point $x\in M$ has a neighborhood that is either flat or
it is a product of a pseudo-Riemannian manifold of constant
sectional curvature $K$ and signature $(r_1,s_1)$ and a
pseudo-Riemannian manifold of constant sectional curvature $-K$
($K\neq 0$) and signature $(r-r_1,s-s_1)$;
\item[4)] ${\rm Hol}^0(M,g)={\rm SO}(r-1,s)$ (resp., ${\rm Hol}^0(M,g) ={\rm SO}(r,s-1)$)
 and each point $x\in M$ has a neighborhood that is either
 flat or
it is a product of a pseudo-Riemannian manifold of constant
sectional curvature and signature $(r-1,s)$ (resp., $(r,s-1)$) and
the space $(L,-(dt)^2)$ (resp., $(L,(dt)^2)$), $L$ is the straight
line or a segment.
\end{itemize}
\end{theorem}

{\bf Proof.} Let $(M,g)$ be a pseudo-Riemannian manifold of
signature $(r,s)$ and dimension $d=r+s$. The vector bundle
$\mathfrak{so}(TM)$ of skew-symmetric endomorphisms of the tangent
bundle $TM$ can be identified with the space of bivectors
$\wedge^2TM$ in such a way that $$(X\wedge Y)Z= g (X,Z)Y- g
(Y,Z)X$$ for all vector fields $X,Y,Z$ on $M$. The Weyl tensor $W$
of the pseudo-Riemannian manifold $(M,g)$ is defined by the
equality
\begin{equation}\label{decR1} W=R+R_L,\end{equation}
where the tensor $R_L$ is defined by
\begin{equation}\label{RL}R_L(X,Y)=LX\wedge Y+X\wedge LY,\end{equation}
$$L=\frac{1}{d-2}\left(\Ric-\frac{s}{2(d-1)}\Id\right)$$ is
the Schouten tensor and $s$ is the scalar curvature.

Suppose that the restricted holonomy group ${\rm Hol}^0(M,g)$ is
not weakly irreducible. The Wu decomposition Theorem \cite{Wu}
states that each point of $M$ has a neighborhood $U$ such that
$(U,g|_U)$ is a product
$$(U,g|_U)=(M_1\times M_2,g_1+g_2)$$
of two pseudo-Riemannian manifolds $(M_1,g_1)$ and $(M_2,g_2)$.
Let $d_1$ and $d_2$ be the dimensions of these manifolds. For the
curvature tensors, Ricci operators and the scalar curvatures it
holds $$R=R_1+R_2,\quad \Ric=\Ric_1+\Ric_2,\quad s=s_1+s_2.$$
First suppose that $d\geq 4$. In this case $W=0$ and we get
\begin{equation}\label{R1R2RL} R_1+R_2=-R_L.\end{equation} Assume that $d_1\geq d_2$ and
$d_1\geq 2$. The curvature tensor $R_1$ can be written in the form
$R_1=W_1-R_{L_1}$. Considering \eqref{R1R2RL} restricted to
$TM_1$, we get that $W_1=0$ and \begin{equation}\label{eqRic12}
\frac{1}{d_1-2}\left(\Ric_1-\frac{s_1}{2(d_1-1)}\id\right)=
\frac{1}{d-2}\left(\Ric_1-\frac{s_1+s_2}{2(d-1)}\id\right).\end{equation}
If $d_2\geq 2$, then taking the trace in \eqref{eqRic12}, we get
$$\frac{s_1}{d_1(d_1-1)}=-\frac{s_2}{d_2(d_2-1)}.$$
 Since $s_1$ is a
function on $M_1$ and $s_2$ is a function on $M_2$, the both
functions must be constant. Substituting the last equality back to
\eqref{eqRic12}, we obtain
\begin{equation}\label{eqRic1}\Ric_1=\frac{s_1}{d_1}\id.\end{equation} Next,
\begin{equation}\label{eqqR1}R_1(X,Y)=\frac{s_1}{d_1(d_1-1)}X\wedge Y.\end{equation} The
same holds for the second manifold. For the sectional curvatures
we get
$$k_1=\frac{s_1}{d_1(d_1-1)}=-\frac{s_2}{d_2(d_2-1)}=-k_2.$$
If $d_2=1$, then \eqref{eqRic12} is equivalent to \eqref{eqRic1}
and this implies \eqref{eqqR1}. From this and the Schur Theorem it
follows that $k_1$ is constant.
 If $d_1=2$, then the curvature
tensor $R_1$ satisfies $R_1(X,Y)=fX\wedge Y$ for some function $f$
on $M_1$. The proof in this case is the same.

If $d=3$, then $d_1=2$ and $d_2=1$. It holds $R=R_1$ and
$R_1(X,Y)=fX\wedge Y$ for some function $f$ on $M_1$. In this case
$(M,g)$ is conformally flat if and only if the Cotton tensor $C$
defined by $$C(X,Y,Z)=g((\nabla_ZL)X,Y)-g((\nabla_YL)X,Z)$$ is
zero. This implies that $f$ is constant, i.e. $(M_1,g_1)$ has
constant sectional curvature.

Now we have to prove that if ${\rm Hol}^0(M,g)$ is irreducible,
then it coincides with ${\rm SO}(r,s)$. Suppose that ${\rm
Hol}^0(M,g)$ is irreducible and it is different from ${\rm
SO}(r,s)$ and ${\rm U}(\frac{r}{2},\frac{s}{2})$. Then the
manifold is Einstein \cite{Bryant2}. Since $(M,g)$ is in addition
conformally flat, $(M,g)$ has constant sectional curvature and its
connected holonomy group must be either trivial or ${\rm
SO}(r,s)$, i.e. we get a contradiction. It is known  that if a
pseudo-K\"ahlerian manifold is conformally flat, then it is flat
\cite{Yano}, hence ${\rm Hol}^0(M,g)\neq {\rm
U}(\frac{r}{2},\frac{s}{2})$. This proves Theorem \ref{ThC1}.
$\Box$

\section{The curvature tensor and the Weyl curvature tensor of Walker
metrics}\label{secCWeyl} In order to prove Theorem \ref{ThC2}, we
give some information about the curvature tensor of the Walker
metric \eqref{Walker}. For the fixed coordinates $v,x^1,...,x^n,u$
consider the field of frames $$p=\partial_v,\quad
X_i=\partial_i-A_i\partial_v,\quad
q=\partial_u-\frac{1}{2}H\partial_v.$$ Consider the distribution
$E$ generated by the vector fields $X_1,...,X_n$. The fibers of
this distribution can be identified with the tangent spaces to the
Riemannian manifolds with the Riemannian metrics $h(u)$. Denote by
$R_0$ the tensor corresponding to  the family of the curvature
tensors of  $h(u)$ under this identification. Similarly denote by
$\Ric(h)$ the corresponding Ricci endomorphism acting on sections
of $E$.

From the results of
\cite{Gal1} it follows that the curvature tensor $R$ of the metric
$g$ can be written in the form \begin{align}\label{eqA} R(p,q)=&
-\lambda p\wedge q-p\wedge \vec{v},\qquad R(X,Y)=R_0(X,Y)-p\wedge
(P(Y)X-P(X)Y),\\ \label{eqB}R(X,q)=&-g(\vec{v},X)p\wedge
q+P(X)-p\wedge T(X),\qquad R(p,X)=0
\end{align}  for all $X,Y\in \Gamma(E)$. Here $\lambda$ is a function, $\vec{v}\in
\Gamma(E)$, $T\in \Gamma({\rm End}(E))$ is symmetric, $T^*=T$,
$R_0=R(h)$, and the tensor $P\in \Gamma(E^*\otimes
\mathfrak{so}(E))$ satisfies
$$g(P(X)Y,Z)+g(P(Y)Z,X)+g(P(Z)X,Y)=0\text{ for all }
X,Y,Z\in \Gamma(E).$$ These elements may be found in terms of the
coefficients of the metric \eqref{Walker}. Let
$P(X_k)X_j=P_{jk}^iX_i$ and $T(X_j)=\sum_iT_{ij}X_j$. Then
$$h_{il}P_{jk}^l=g(R(X_k,q)X_j,X_i),\quad T_{ij}=-g(R(X_i,q)q,X_j).$$
Using direct computations, we obtain
\begin{align}\label{lamv}\lambda=&\frac{1}{2}\partial^2_vH,\quad
\vec{v}=\frac{1}{2}\left(\partial_i\partial_vH-A_i\partial^2_vH\right)h^{ij}X_j,\\
\label{eqPijk}h_{il}P_{jk}^l=&-\frac{1}{2}\nabla_kF_{ij}+\frac{1}{2}\nabla_k\dot{h}_{ij}-\dot\Gamma^l_{kj}h_{li},\\
\label{eqTij}
T_{ij}=&\frac{1}{2}\nabla_i\nabla_jH-\frac{1}{4}(F_{ik}+\dot
h_{ik})(F_{jl}+\dot
h_{jl})h^{kl}-\frac{1}{4}(\partial_vH)(\nabla_iA_j+\nabla_jA_i)
\\ \nonumber &-\frac{1}{2}(A_i\partial_j\partial_vH+A_j\partial_i\partial_vH)-\frac{1}{2}(\nabla_i\dot A_j+\nabla_j\dot A_i)\\
\nonumber
&+\frac{1}{2}A_iA_j\partial_v^2H+\frac{1}{2}\textit{\"{h}}_{ij}+\frac{1}{4}\dot
h_{ij}\partial^2_vH,
\end{align} where $$F=dA,\quad F_{ij}=\partial_iA_j-\partial_jA_i$$
is the differential of the 1-form $A$, and the covariant
derivatives are taken with  respect to the metric $h$. In the case
of $h$, $A$ and $H$ independent of $u$, the curvature tensor of
the metric \eqref{Walker} is found in \cite{G-P}.

The Ricci operator has the following form:
\begin{align}\label{Ric3} \Ric(p)=&\lambda p,\quad \Ric(X)=-g(X,\tRic P-\vec{v})p+\Ric(h)(X),\\
\label{Ric4} \Ric(q)=&-(\tr T)p-\tRic(P)+\vec{v}+\lambda
q,\end{align} where $\tRic P=h^{ij}P(X_i)X_j$ \cite{onecomp}. For
the scalar curvature we get $s=2\lambda+s_0$, where $s_0$ is the
scalar curvature of $h$. Using this, we may compute the tensor
$R_L$,
\begin{align}\label{eq12}
R_L(p,X)&=\frac{1}{n}p\wedge\left(\Ric(h)+\frac{(n-1)\lambda-s_0}{n+1}\id\right)
X,\\\label{eq13}
R_L(p,q)&=\frac{1}{n}\left(\frac{2n\lambda-s_0}{n+1}p\wedge
q+p\wedge (\vec{v}-\tRic P)\right),\\ \label{eq14}
R_L(X,Y)&=\frac{1}{n}\left(p\wedge(g(X,\vec{v}-\tRic P)Y-g(Y,\vec{v}-\tRic
P)X)\right.\\ \nonumber
&\left.+\left(\Ric(h)-\frac{s}{2(n+1)}\right)X\wedge Y+X\wedge
\left(\Ric(h)-\frac{s}{2(n+1)}\right)Y\right),\\
\label{eq15} R_L(X,q)&=\frac{1}{n}\Big((\tr T)p\wedge X+
g(X,\vec{v}-\tRic P)p\wedge q+X\wedge(\vec{v}-\tRic
P)\\\nonumber&+\left(\Ric(h)+\frac{(n-1)\lambda-s_0}{n+1}\id\right)X\wedge
q\Big).
\end{align}
The Weyl tensor $W$ can be computed using this and \eqref{decR1}.

\section{Proof of Theorem \ref{ThC2}}\label{secCdokvo}
Suppose that $(M,g)$ is a Walker manifold, i.e. its holonomy
algebra is contained in $\mathfrak{sim}(n)$. The local form of the
metric $g$ is given by
 \eqref{Walker}. Suppose that $g$ is conformally flat, i.e.
$W=0$.

\begin{lemma}\label{lemC1} The equation $W=0$ is equivalent to the following system of equations:
\begin{equation}\label{condW0} s_0=-n(n-1)\lambda,
\quad R_0=-\frac{1}{2}\lambda R_{\rm id},\quad
P(X)=\vec{v}\wedge X,\quad T=f\id_E,\end{equation} where $X$ is
any section of $E$ and $f$ is a function. In particular, $W=0$ implies that
$\tRic P=-(n-1)\vec{v}$ and the Weyl tensor $W_0$ of $h$ is zero.
\end{lemma}

{\bf Proof.} Suppose that $W=0$. Then it holds $R=-R_L$. From \eqref{eqB} and \eqref{eq12}
it follows that $\Ric(h)=-\frac{(n-1)\lambda-s_0}{n+1}\id$. Taking
the trace, we get  $s_0=-n(n-1)\lambda$. Hence,  $\Ric(h)=\frac{s_0}{n}\id$, i.e. the metrics $h$ are Einstein,
and it holds
$$R_0=W_0-\frac{s_0}{2n(n-1)}R_{\rm id}.$$
 From \eqref{eqA} and
\eqref{eq14} it follows that
$$R_0+\frac{s_0}{2n(n-1)}R_{\rm id}=0.$$
We conclude that
$W_0=0$. Using \eqref{eqB} and
\eqref{eq15}, we get $P(X)=-\frac{1}{n}X\wedge (\vec{v}-\tRic P)$.
Applying $\tRic$, we obtain $\tRic P=-(n-1)\vec{v}$. Consequently,
$P(X)=\vec{v}\wedge X$. From \eqref{eqB} and \eqref{eq15} it
follows that $T(X)=\frac{1}{n}(\tr T)X$, i.e. $T=f\id$ for a
function $f$. Conversely, \eqref{condW0} implies $W=0$. $\Box$

 Since $h$ is independent of $v$, $\partial_v\lambda=0$.
From \eqref{lamv} it follows that $$H=\lambda v^2+H_1v+H_0,\quad
\partial_vH_1=\partial_v H_0=0.$$ From \eqref{eqPijk} it follows that
the components of tensor $P$    do not depend on the coordinate
$v$. This,  the equation $P(X)=\vec{v}\wedge X$ and \eqref{lamv}
imply that $\partial_i\partial^2_v H=0$, i.e.
$\partial_i\lambda=0$, consequently $s_0$ and $\lambda$ are
functions depending only on $u$. Note that for $n\geq 3$ this
follows also from $\Ric(h)=\frac{s_0}{2}\id$. We conclude that
each metric in the family $h(u)$ is of constant sectional
curvature. There exists a transformation $\tilde x^i=\tilde
x^i(x^k,u)$ of the coordinates $x^1,...,x^n$ depending on the
parameter $u$ such that with respect to the new coordinates the
metric $h$ takes the form
$$h=\Psi\sum_{k=1}^n(dx^k)^2,\quad \Psi=\frac{4}{\left(1-\lambda(u)\sum_{k=1}^n(x^k)^2\right)^2}.$$
Considering the transformation $\tilde v=v,\tilde x^i=\tilde x^i(x^k,u),\tilde u=u$, we may assume that $h$ in \eqref{Walker}
has the just obtained form.

Let us consider the equation $T_{ij}=f\delta_{ij}$. From
\eqref{eqTij} it follows that $f=vf_1+f_0$,
$\partial_vf_1=\partial_v f=0$. Applying $\partial_v$ to
$T_{ij}=f\delta_{ij}$, we get the equations
$$f_1\delta_{ij}=\frac{1}{2}\nabla_i\nabla_jH_1-\lambda(u)\frac{1}{2}(\nabla_iA_j+\nabla_jA_i).$$
These equations are equivalent to the equations \begin{equation}\label{cfeaonA}\nabla_iZ_i=\nabla_jZ_j,\quad \nabla_iZ_j+\nabla_jZ_i=0,\qquad i\neq j,\end{equation}
where \begin{equation}\label{eqdZ} Z_i=\lambda A_i-\frac{1}{2}\partial_iH_1.\end{equation}
Note that $\vec v=-Z_ih^{ij}X_j$.
The Christoffel symbols of the metric $h$ are the following:
$$\Gamma^k_{ij}=\frac{1}{2\Psi}(\delta_{kj}\partial_i\Psi+\delta_{ki}\partial_j\Psi-\delta_{ij}\partial_k\Psi).$$
Using that, the above equations may be rewritten in the form
\begin{equation}\label{cf3} \partial_i\left(\frac{Z_i}{\Psi}\right)=\partial_j\left(\frac{Z_j}{\Psi}\right),\quad
\partial_i\left(\frac{Z_j}{\Psi}\right)+\partial_j\left(\frac{Z_i}{\Psi}\right)=0,\quad i\neq j.\end{equation}

We will distinguish the cases $n=2$ and $n\geq 3$.

{\bf Case  $n\geq 3$.}

\begin{lemma}\label{cfLemsol}
 If $n\geq 3$, then the general solution of the system
$$\partial_if_i=\partial_jf_j,\quad \partial_if_j+\partial_jf_i=0,\quad i\neq j$$
has the form
$$f_i=x^iB_kx^k-\frac{1}{2}B_i\sum_{k=1}^n(x^k)^2+d_{ik}x^k+cx^i+c_i,$$
where $B_i,c_i,c,d_{ik}\in\mathbb{R}$, $d_{ki}=-d_{ik}$.
\end{lemma}

{\it Proof.} Let $i,j,k$ be pairwise different, then
$$\partial_i\partial_jf_k=-\partial_i\partial_k f_j=\partial_k\partial_jf_i=-\partial_i\partial_jf_k,$$
i.e. $\partial_i\partial_j f_k=0$. This shows that
$f_k=\sum_{i\neq k}C_{ki}(x^i,x^k)$. Then it is not hard to find
these functions. $\Box$

We conclude that
\begin{equation}\label{cfSAi=}
Z_i=\Psi \left(x^iB_k(u)x^k-\frac{1}{2}B_i(u)\sum_{k=1}^n(x^k)^2+d_{ik}(u)x^k+c(u)x^i+c_i(u)\right),
\end{equation}
where $B_i(u),c_i(u),c(u),d_{ik}(u)$ are functions of $u$, and
$d_{ki}(u)=-d_{ik}(u)$.

 The  system of equations that we have solved is very similar to the equation of the Killing 1-form:
 $$\nabla_i Z_j+\nabla_j Z_i=0.$$
Let us prove the following lemma.

\begin{lemma}\label{LemKilS} Any Killing vector field on the space with the metric $\Psi\sum_{k=1}^n(dx^k)^2$ has the following form $$X=X^i\partial_i,\quad X^i=f_{ik}x^k -2\lambda x^ic_kx^k+\lambda c_i\sum_{k=1}^n(x^k)^2+c_i,$$
where $b_i,f_{ik}\in\mathbb{R}$, $f_{ki}=-f_{ik}$. \end{lemma}

Depending on the value of $\lambda$, the above metric is the metric of one of the spaces: the sphere, the Lobachevskian space, the Euclidean space. This description corresponds to the fact that the Lie algebra $\mathfrak{k}$
 of the Killing vector fields on these
spaces  is isomorphic respectively  to $\mathfrak{so}(n+1)$,
$\mathfrak{so}(1,n)$, $\mathfrak{iso}(\mathbb{R}^n)$.
 The symmetric decomposition of the Lie algebra $\mathfrak{k}$ is of the form $\mathfrak{k}=\mathfrak{so}(n)+\mathbb{R}^n$. The vector fields defined by the numbers $f_{ik}$ correspond to elements of $\mathfrak{so}(n)$, while the vector fields  defined by the numbers $c_i$   correspond to elements of $\mathbb{R}^n$.

 {\it Proof of Lemma \ref{LemKilS}.}
 Consider the Killing 1-form $Z_i=h_{ij}X^j$. In addition to equations \eqref{cfeaonA}  it satisfies the equations $\nabla_iZ_i=0$. These equations take the form
 $$\partial_i\left(\frac{Z_i}{\Psi}\right)=-\frac{1}{2\Psi^2}\sum_{k=1}^nZ_k\partial_k\Psi.$$
 This implies that $Z_i=\Psi f_i$, where $f_i$ are given by Lemma \ref{cfLemsol}  with $c=0$ and $B_k=-2\lambda c_k$.
 This proves the Lemma. $\Box$

  From \eqref{eqdZ} it follows that
   \begin{equation}\label{lFij} \lambda F_{ij}=\partial_iZ_j-\partial_jZ_i=\Psi^{\frac{3}{2}}\big((B_i-2\lambda C_i)x^j-(B_j-2\lambda C_j)x^i
+2\lambda x^k(d_{jk}x^i-d_{ik}x^j)\big)-2\Psi d_{ij}.
   \end{equation}
   This easily implies that $B_i(u)=\lambda(u)\tilde B_i(u)$ for some functions $\tilde B_i(u)$.

From \eqref{eqPijk} and \eqref{lamv} it follows that
the equation $P(X)=\vec{v}\wedge X$ takes the form
\begin{multline}
\label{eqPvX}
-\frac{1}{2}\nabla_kF_{ij}-\delta_{ik}\frac{\Psi}{2}\partial_u\partial_j\ln\Psi+\delta_{jk}\frac{\Psi}{2}\partial_u\partial_j\ln\Psi\\=
\Psi \left(\frac{1}{2}\partial_jH_1-\lambda
A_j\right)\delta_{ki}-\Psi \left(\frac{1}{2}\partial_iH_1-\lambda
A_i\right)\delta_{kj}.\end{multline} If $i$, $j$, $k$ are
pair-wise different, then
\begin{multline*}0=\nabla_kF_{ij}=\partial_kF_{ij}-\frac{1}{\Psi}F_{ij}\partial_k\Psi-\frac{1}{2\Psi}F_{kj}\partial_i\Psi-
\frac{1}{2\Psi}F_{ik}\partial_j\Psi\\=\Psi\partial_k\left(\frac{F_{ij}}{\Psi}\right)-\lambda
F_{kj}x^i\sqrt{\Psi}-\lambda F_{ik}x^j\sqrt{\Psi}.\end{multline*}
Using this and Equation \eqref{lFij}, we obtain
$$-\partial_k\left(\frac{F_{ij}}{\Psi}\right)=\lambda x^k\Psi\big((\tilde B_j-2C_j)x^i-(\tilde B_i-2C_i)x^j\big)+
2\sqrt{\Psi}(x^id_{kj}-x^jd_{ki})\big)+2\lambda\Psi x^k(x^jd_{il}-x^id_{jl})x^l.$$
Integrating over $x^k$, we get
\begin{equation}\label{eqFij}
F_{ij}=\Psi^{\frac{3}{2}}\big((\tilde B_i-2C_i)x^j-(\tilde B_j-2C_j)x^i+2(d_{li}x^j-d_{lj}x^i)x^l\big)-\Psi C_{ij}(x^i,x^j,u)
\end{equation}   for some functions $C_{ij}(x^i,x^j,u)$.
Comparing this with \eqref{lFij}, we conclude that $d_{ij}=\frac{\lambda}{2}C_{ij}$.
Equation \eqref{eqPvX} for $k=i\neq j$ is of the form
$$\nabla_iF_{ij}=-\Psi\partial_u\partial_j\ln\Psi+2\Psi Z_j.$$
Direct computations show that
$$\nabla_iF_{ij}=\Psi^{\frac{3}{2}}\partial_i\left(\frac{F_{ij}}{\Psi^{\frac{3}{2}}}\right)+\lambda F_{lj}x^l\sqrt{\Psi}.$$
Using this, \eqref{eqdZ}, \eqref{lFij} and \eqref{eqFij}, we get
$$\Psi^{\frac{3}{2}}\left(-2d_{ij}-\partial_i\frac{C_{ij}}{\sqrt\Psi}\right)+\Psi\partial_u\partial_j\ln\Psi-\Psi^2\left(
\lambda D_kx^kx^j-\frac{1}{2}D_j\sum_{k=1}^n(x^k)^2+2c(u)x^j+\frac{1}{2}D_j\right)=0,$$
where $D_k=\tilde B_k+2C_k$. Using the equalities
$$\partial_j\partial_u\ln\Psi=\dot\lambda x^j\Psi,\quad   \partial_i\frac{C_{ij}}{\sqrt\Psi}=-\lambda x_iC_{ij}+\frac{\partial_iC_{ij}}{\sqrt \Psi},\quad
d_{ij}=\frac{\lambda}{2}C_{ij},$$ we get
$$-\partial_iC_{ij}-\Psi\left(
\lambda
D_kx^kx^j-\frac{1}{2}D_j\sum_{k=1}^n(x^k)^2+(2c(u)-\dot\lambda)x^j+\frac{1}{2}D_j\right)=0.$$
If $\lambda(u)\neq 0$, then the equality
$d_{ij}=\frac{\lambda}{2}C_{ij}$ implies $\partial_iC_{ij}=0$.
Consequently, $D_k=0$ and $c=\frac{1}{2}\dot\lambda$. If
$\lambda(u)=0$, then we get
$$\partial_j\partial_iC_{ij}+4(D_jx_j-(2c-\dot\lambda))=0.$$ Here we wrote $x_j$ instead of $x^j$
in order to avoid the sum over $j$. Since $C_{ij}=-C_{ji}$, we get
$D_j=0$, $c=\frac{1}{2}\dot\lambda$, and $\partial_iC_{ij}=0$.
Thus it holds $D_k=0$, $\tilde B_k=-2C_k$,
$c=\frac{1}{2}\dot\lambda$, and $d_{ij}=\frac{\lambda}{2}C_{ij}$,
where $C_{ij}$ are functions of $u$. We conclude that
\begin{equation}
\label{Fij1}
 F_{ij}=\Psi^{\frac{3}{2}}\big(4C_j(u)x^i-4C_i(u)x^j+\lambda(u)(C_{li}(u)x^j-C_{lj}(u)x^i)x^l\big)-\Psi C_{ij}(u).\end{equation}
Recall that $F=dA$. In \cite{G-P} it is noted that the
transformation $v\mapsto v-\phi(x^1,...,x^m,u)$ changes $A_i$ to
$A_i+\partial_i\phi$. Clearly, this transformation does not change
$F$. This allows us to choose any $A$ such that $dA=F$. We take
$$A_i=\Psi\left(-4C_k(u)x^kx^i+2C_i(u)\sum_{k=1}^n(x^k)^2+\frac{1}{2}C_{ik}(u)x^k\right).$$
Consider the coordinate transformation with the inverse one
\begin{equation}\label{cf1}
v=\tilde v,\quad x^i=A^i_j(\tilde u)\tilde x^j,\quad u=\tilde u,
\end{equation}
where $A^i_j(u)$ is a family of orthogonal matrices. It is easy to
check that
$$\tilde A_i=\sum_{k=1}^n A^k_i(u)(\partial_{u}A^k_l(u))\tilde x^l+A^k_i(u)A_i.$$
The obtained metric has the same form and it holds  $$\tilde
C_i(u)=C_j(u)A^j_i(u),\quad \tilde C_{ij}(u)=\sum_{k=1}^n
A^k_i(u)\partial_{
u}A^k_j(u)+\frac{1}{2}A^r_i(u)C_{rk}(u)A^k_j(u).$$ Consider the
equation $\tilde C_{ij}(u)=0$. Since $\sum_{k=1}^n
A_k^i(u)A_k^j(u)=\delta^{ij}$, it can be written in the form
$$\partial_uA^k_i(u)=A^j_i(u) \frac{1}{2}C_{jk}(u).$$ Since $C_{jk}(u)$
is skew-symmetric, $\frac{1}{2}C_{jk}(u)$ is a curve in the Lie
algebra $\mathfrak{so}(n)$. Then $A^k_i(u)$ satisfying the above
equation is nothing else as the development of the curve
$\frac{1}{2}C_{jk}(u)$ in the Lie group ${\rm SO}(n)$. Thus,
applying such transformation, we may assume that $C_{ij}(u)=0$.

Next, $$\partial_iH_1=2\lambda A_i-2Z_i=-2\Psi\left(2\lambda
C_k(u)x^kx^i- \lambda
C_i(u)\sum_{k=1}^n(x^k)^2+\frac{1}{2}\dot\lambda
x^i+C_i(u)\right).$$ We conclude that
$$H_1=-4C_k(u)x^k\sqrt{\Psi}-\partial_u\ln\Psi+K(u)$$ for some function
$K(u)$.

We are left with the only unknown function $H_0$. Consider
equations $T_{ij}=f\delta_{ij}$. If $i\neq j$, then
$\nabla_i\nabla_jH_0=\sqrt{\Psi}\partial_i\partial_j\frac{H_0}{\sqrt\Psi}$,
and using \eqref{eqTij}, we obtain
$$\partial_i\partial_j\frac{H_0}{\sqrt\Psi}=2\Psi^{\frac{3}{2}}x^ix^j\sum_{k=1}^nC_k^2(u).$$
If $\lambda(u)\neq 0$, then the function
$\frac{4}{\lambda^2(u)}\Psi\sum_{k=1}^nC_k^2(u)$ is a partial
solution of the system. The general solution is of the form
$$H_0=\frac{4}{\lambda^2(u)}\Psi\sum_{k=1}^nC_k^2(u)+\sqrt\Psi\sum_{k=1}^nf_k(x^k,u).$$
The condition $\nabla_i\nabla_iH_0=\nabla_j\nabla_jH_0$  implies
$\partial_i^2f_i=\partial_j^2f_j$, hence
$$f_i(x^i,u)=a(u)(x^i)^2+D_i(u)x_i+d_i(u)$$ for some functions
$a(u)$, $D_i(u)$, $d_i(u)$.  We see that $H_0$ is as in the
statement of the theorem for the case $\lambda(u)\neq 0$. The case
$\lambda(u)=0$ is similar.

{\bf Case $n=2$.}  The system of equations \eqref{cf3} takes the form
\begin{equation}\label{cf4} \partial_1\left(\frac{Z_1}{\Psi}\right)=\partial_2\left(\frac{Z_2}{\Psi}\right),\quad
\partial_1\left(\frac{Z_2}{\Psi}\right)+\partial_2\left(\frac{Z_1}{\Psi}\right)=0\end{equation}
that implies  only that $\frac{Z_1}{\Psi}$ and $\frac{Z_2}{\Psi}$ are real and imaginary parts of a complex homomorphic function of the variable $x^1+ix^2$.

Note  that
$$\nabla_1F_{12}=\Psi\partial_1\frac{F_{12}}{\Psi},\quad
\nabla_2F_{12}=\Psi\partial_2\frac{F_{12}}{\Psi}.$$ Using
\eqref{eqPvX}, we get
$$Z_1=-\frac{1}{2}\partial_2\frac{F_{12}}{\Psi}+\frac{1}{2}\partial_u\partial_1\ln\Psi,\quad Z_2=\frac{1}{2}\partial_1\frac{F_{12}}{\Psi}+\frac{1}{2}\partial_u\partial_2\ln\Psi.$$
Substituting that to the first equation in \eqref{cf4}, and using
the equality $\partial_u\partial_i\ln\Psi=\dot\lambda x^i\Psi$, we
get
$$\partial_1\partial_2\frac{F_{12}}{\Psi^{\frac{3}{2}}}=0.$$ This implies
$$F_{12}=\Psi^{\frac{3}{2}}(f_1(x^1,u)+f_2(x^2,u)).$$ Substituting
that to the second equation in \eqref{cf4}, we obtain
$\partial_1^2f_1=\partial_2^2f_2$, i.e.
$$f_1=t(u)(x^1)^2+a^1(u)x^1+b_1(u),\quad f_2=t(u)(x^2)^2+a^2(u)x^2+b_2(u).$$
It holds $$\partial_1H_1=2Z_1+2\lambda A_1,\quad
\partial_2H_1=2Z_2+2\lambda A_2.$$ Consequently,
$$0=\partial_2\partial_1H_1-\partial_1\partial_2H_1=-2\lambda F_{12}+2\partial_2Z_1-2\partial_1Z_2.$$
Direct computations show that this implies $t=\lambda b$. Thus,
$$F_{12}=\Psi^{\frac{3}{2}}(\lambda(u)b(u)((x^1)^2+(x^2)^2)+a_i(u)x^i(u)+b(u)).$$
Equality \eqref{Fij1} for $n=2$ can be rewritten in the form
$$F_{12}=\Psi^{\frac{3}{2}}\left(-\frac{1}{2}\lambda(u)C_{12}(u)((x^1)^2+(x^2)^2)-\frac{1}{2}C_{12}(u)+4C_2(u)x^1-4C_1(u)x^2\right).$$ This shows that $F_{ij}$ is the same as for $n\geq 3$, and
 the rest of the proof for $n=2$ is the same as for $n\geq 3$.
 The theorem is true.
$\Box$

\section{Proof of Theorem \ref{ThC2A}}\label{secpT2A}
Consider the coordinates $v,x^1,...,x^n,u$ and the metric $g$ is in Theorem \ref{ThC2}.

{\bf 1)} Suppose that $\lambda(u)$ is nowhere vanishing. Consider
the transformation $$v\mapsto v-\phi,\quad
\phi=\frac{2}{\lambda(u)}C_k(u)x^k\sqrt\Psi.$$ Then $A_i$ changes
to
$$A_i+\partial_i\phi=\Psi\left(-2C_k(u)x^kx^i+C_i(u)\sum_{k=1}^n(x^k)^2+\frac{C_i(u)}{\lambda(u)}\right).$$
By Lemma \ref{LemKilS}, for each $u$, $h^{ik}A_k$ is a Killing
vector field of the Riemannian metric $h(u)$.

Following the ideas from \cite{GL10}, we are looking for a
coordinate transformation in order to set $A_i$ to zero. Consider
the coordinate transformation with the inverse one given by
$$v=\tilde v,\,x^k=x^k(\tilde x^i,u),\,u=\tilde u.$$ It holds
$$\tilde A_k=\frac{\partial x^i}{\partial\tilde
x^k}\left(A_i+h_{ij}\frac{\partial x^j}{\partial u}\right).$$ The
equation $\tilde A_k=0$ is of the form
$$ \frac{\partial x^j(\tilde x^1,...,\tilde x^n,u)}{\partial u}=V^j(x^1(\tilde x^1,...,\tilde x^n,u),...,x^n(\tilde x^1,...,\tilde x^n,u),u),$$ where $V^j=-h^{jk}A_k$. Considering $\tilde x^1,...,\tilde x^k$ as parameters and imposing the initial dates
$$x^i(\tilde x^1,...,\tilde x^n,u_0)=\tilde x^i,$$ we see that the above system is a system of ordinary differential equations depending on the initial dates as on parameters. Such system has a unique solution that gives the required transformation.
Since for each $u$, $h^{ik}A_k$ is a Killing vector field of the
Riemannian metric $h(u)$, the transformation $x^k=x^k(\tilde
x^i,u)$ is an isometry, hence $h$ remains the same. Thus we get
the metric from Theorem \ref{ThC2} with $C_i(u)=0$. Applying the
transformation $v\mapsto v+\frac{1}{2\lambda(u)}K(u)$, we get
$H_1=-\partial_u\ln\Psi$.

{\bf 2)} Suppose that $\lambda$ is constantly zero. The from $F$  is given by
$$F_{ij}=32(C_j(u)x^i-C_i(u)x^j).$$
As it is explained  in the proof of Theorem \ref{ThC2}, we may
take $A_i=16C_i(u)\sum_{k=1}^n(x^k)^2$. Considering the
transformation $x^i\mapsto \frac{x^i}{2}$, and redenoting
$2C_i(u)$ by $C_i(u)$, we get $A=C_i(u)\sum_{k=1}^n(x^k)^2dx^idu$,
i.e. $A_i=C_i(u)\sum_{k=1}^n(x^k)^2$. Moreover,
$h=\sum_{k=1}^n(dx^k)^2$. From the equation $P(X)=\vec v\wedge X$
it follows that $$\partial_jH_1=-\nabla_iF_{ij}=-2C_j.$$ Hence,
$H_1=-2C_kx^k+K(u)$ for some function $K(u)$.

Suppose that $\sum_{k=1}^n C^2_k\not\equiv 0$. Consider the
coordinate transformation with the inverse one
$$v=\tilde v,\quad x^i=\tilde x^i+b^i(\tilde u),\quad u=\tilde u$$
such that $2C_i(u)b^i(u)=K(u)$. After that $H_1=-2C_i(u)x^i$, i.e.
we may assume that $K(u)=0$. The equation $T_{ij}=f\delta_{ij}$ for $i\neq j$ takes the form
$$\frac{1}{2}\partial_i\partial_jH_0=x^ix^j\sum_{k=1}^nC_k-2C_kx^k(C_jx^i+C_ix^j)-C_iC_j\sum_{k=1}^n(x^k)^2+(\dot C_jx^i+\dot C_ix^j). $$
The general solution to this system under the condition
$T_{ii}=T_{jj}$ is given in the statement of the theorem.

Suppose that $\sum_{k=1}^n C^2_k\equiv0$. It is easy to check that
the non-zero components of the curvature tensor are $R_{iuju}$,
$R_{uiuj}$, $R_{iuuj}$, $R_{uiju}$, hence  the metric is a pp-wave
\cite{ESI}, i.e. $h=\sum_{k=1}^n(dx^k)^2$, $A=0$, and $H=H_0$.
Moreover, $$H_0=a(u)\sum_{i=1}^n(x^i)^2+D_i(u)x^i+D(u).$$ Consider
the new coordinates $$\tilde v = v  -\sum_j\frac{d b^j(u)}{d
u}x^j+d(u), \quad \tilde x^i = x^i+b^i(u), \quad \tilde u = u.$$
We obtain the metric of the same form with
$$\tilde H_0=a(u)\sum_{i=1}^n(x^i)^2+\tilde D_i(u)x^i+\tilde D(u),$$
where \begin{align} \label{tGk}\tilde D_j&=-2\frac{d^2 b^j}{(d
 u)^2}+2ab^j+D_j,\\
\tilde D&=2\frac{dd( u)}{d u}
+\sum_j\left(\frac{db^j}{du}\right)^2+a\sum_j(b^j)^2+D_ib^i+D.\end{align}
Equation \eqref{tGk} implies the existence of  $b^j(u)$ such that
$\tilde D_k=0$. Using the last equation, we can chose $d(u)$ such
that $\tilde D=0$. Thus,
$$H_0=a(u)\sum_{i=1}^n(x^i)^2.$$
The theorem is true. $\Box$

\section{Proof of Theorem \ref{ThC3}}\label{secpT3}
We consider the  metric $g$ from Theorem \ref{ThC2} and consider different cases.

Suppose that $\lambda\equiv 0$ on $M$. Suppose that $\sum_{k=1}^nC_k^2\not\equiv 0$. Then there exists a point
$x\in M$  such that $\vec{v}_x\neq 0$. The condition on
the curvature tensor shows that
$$R_x(p_x,q_x)=-p_x\wedge \vec{v}_x,\quad R_x(X,Y)=p_x\wedge ((X\wedge Y)\vec{v}_x).$$
This shows that $p_x\wedge E_x\subset \mathfrak{g}$. Next,
$$R_x(\vec{v}_x,q_x)=-g(\vec{v}_x,\vec{v}_x)p_x\wedge q_x-p_x\wedge
T_x(\vec{v}_x),$$ which implies $\mathbb{R} p_x\wedge
q_x\subset\mathfrak{g}$. Finally,
$$R_x(X,q_x)=-g(\vec{v}_x,X)p_x\wedge q_x+\vec{v}_x\wedge
X-p_x\wedge T_x(X).$$ Since the bivectors of the form
$\vec{v}_x\wedge X$ generate the Lie algebra $\mathfrak{so}(E_x)$,
we conclude that
$$\mathfrak{g}=\mathbb{R} p_x\wedge q_x+\mathfrak{so}(E_x)+p_x\wedge E_x\simeq\mathfrak{sim}(n).$$
If  $\sum_{k=1}^nC_k^2\equiv 0$ and $a^2\not\equiv 0$, then  the
metric is given by \eqref{cfppw}, and its holonomy algebra
coincides with $\mathbb{R}^n$. If $\sum_{k=1}^nC_k^2(u)\equiv0$
and $a^2(u)\equiv 0$, then the metric is flat.

Suppose that $\lambda$ is a non-zero constant, then the metric can
be written as the first metric from Theorem \ref{ThC2A}. It holds
$$R(p,q)=-\lambda p\wedge q,\quad R(X,Y)=-\frac{1}{2}\lambda
X\wedge Y,\quad X,Y\in\Gamma(E).$$ This shows that $\mathfrak{g}$
contains the subalgebra
$\mathbb{R}\oplus\mathfrak{so}(n)\subset\mathfrak{sim}(n)$. Next,
$$T_{ii}=\lambda
\Psi^{\frac{3}{2}}D_kx^k+\Psi(\sqrt\Psi-1)(a+\lambda D).$$ If
$\sum_{k=1}^nD^2_k+(a+\lambda D)^2\not\equiv 0$, then the tensor
$T$ is not identical zero, and the equality $R(X,q)=-n\tr T
p\wedge X$ shows that $\mathfrak{g}$ contains
$\mathbb{R}^n\subset\mathfrak{sim}(n)$, i.e.
$\mathfrak{g}=\mathfrak{sim}(n)$. Otherwise,
$$g=\Psi\sum_{k=1}^n(dx^k)^n+2dvdu+(\lambda v^2+2D(u))(du)^2.$$ We
see that $g$ is decomposable. The metric $2dvdu+(\lambda
v^2+2D(u))(du)^2$ is of constant sectional curvature
$\frac{\lambda}{2}$, hence it is isometric to the metric
$2dvdu+\lambda v^2(du)^2$.

Suppose that $\dot\lambda\not\equiv 0$ for some coordinate system.
Then $\lambda\not\equiv 0$ on an open subset of the definition
domain of this system, and the metric can be written as the first
metric from Theorem \ref{ThC2A}. It holds $\vec
v=-\frac{1}{2}\dot\lambda x^i\partial_i$. Let $x$ be any point
such that $\vec v_x\neq 0$. It holds $$R_x(p_x,q_x)=-\lambda_x
p_x\wedge q_x-p_x\wedge \vec v_x\in\mathfrak{g},$$
$$R_x(X,Y)=-\frac{\lambda}{2} X\wedge Y+p_x\wedge((Y\wedge X)\vec v_x)\in \mathfrak{g},\quad X,Y\in T_xM.$$
Taking the Lie bracket of the last two elements of $\mathfrak{g}$,
we obtain $p_x\wedge ((Y\wedge X)\vec v_x)\in \mathfrak{g}$. This
shows that $\mathfrak{g}$ contains the subalgebra isomorphic to
$\mathbb{R}^n\subset\mathfrak{g}$. The above two equalities imply
that $\mathfrak{g}$ contains the subalgebra
$\mathbb{R}\oplus\mathfrak{so}(n)\subset\mathfrak{sim}(n)$. Thus,
$\mathfrak{g}=\mathfrak{sim}(n)$.

The theorem is true. $\Box$

\section{The Ricci operator of the obtained
metrics}\label{secCRic}

The Ricci operator of the  metric \eqref{cfppw} has the form
$$\Ric=na(u)\partial_v\otimes du,$$ in particular, $\Ric^2=0$.

In \cite{E-I}, complete conformally flat Lorentzian  manifolds
$(M,g)$ satisfying the condition
\begin{equation}\label{cond1}[R(X,Y),\Ric]=0\end{equation}
are studied. It is shown that these manifolds are exhausted by the
spaces of constant sectional curvature, by the products of two
spaces of constant sectional curvature, and by products of  spaces
of constant sectional curvature with  intervals. The Ricci
operator of the  metric  \eqref{cfppw} satisfies \eqref{cond1}.
Moreover, such metric is complete, e.g. for $a(u)=1$, i.e. for the
Cahen-Wallach spaces. Thus some metrics in
 \cite{E-I} are loosen. In
\cite{Hon}, pseudo-Riemannian conformally flat manifolds $(M,g)$
satisfying \eqref{cond1} are studied. It is shown that in addition
to the obvious cases, $(M,g)$ may be a complex sphere or a space
satisfying $\Ric^2=0$.  Various examples of conformally flat
manifolds with $\Ric^2=0$ are constructed in \cite{H-K}.

The Ricci operator of the second metric from
Theorem \ref{ThC2A} is the following:
$$ \Ric(p)=0,\quad \Ric(X)=ng(X,\vec v)p,\quad
\Ric(q)=nT_{11}p+n\vec v, \quad X\in\Gamma(E),$$ here $\vec
v=-\sum_{k=1}^nC_k(u)X_k$. The function $T_{11}$ can be found
using \eqref{eqTij}. It holds $\Ric^2\neq 0$ and $\Ric^3=0$.
Condition \eqref{cond1} is not satisfied. The scalar curvature of
this metric is zero.

For the first metric form Theorem \ref{ThC2A} it holds
$$\Ric(p)=\lambda p,\quad \Ric(X)=ng(X,\vec{v})p-(n-1)\lambda X,\quad
\Ric(q)=nT_{11}p+n\vec{v}+\lambda q,$$
where $\vec{v}=-\frac{1}{2}\dot\lambda\Psi x^kX_k$.
The Ricci operator is not nilpotent.
The scalar curvature equals to $2\lambda+s_0=-(n-2)(n+1)\lambda$ and it is zero only in dimension four.

\section{The case of dimension 4}\label{secCdim4}

 Applying Theorem \ref{ThC1} to a  conformally flat nonflat
 Lorentzian manifold $(M,g)$ of dimension 4 with the holonomy algebra $\mathfrak{g}\subset\mathfrak{so}(1,3)$,
 we obtain that $(M,g)$ must satisfy one of the following conditions:
\begin{description}
\item[1)] $\mathfrak{g}=\mathfrak{so}(1,3)$;
\item[2)] $\mathfrak{g}\subset\mathfrak{sim}(2)$, i.e. $(M,g)$ is as in Theorem \ref{ThC2} with $n=2$;
\item[3)] $\mathfrak{g}=\mathfrak{so}(1,1)\oplus\mathfrak{so}(2)$, and $(M,g)$ is locally isometric either to the product of
$(dS_2,cg_{dS_2)}$ with $(L^2,cg_{L^2})$, or to the product of
$(AdS_2,cg_{AdS_2})$ with $(S^2,cg_{S^2})$;
\item[4)] $\mathfrak{g}=\mathfrak{so}(1,2)$, and $(M,g)$ is locally isometric to the product of $(\mathbb{R},(dt)^2)$ either with
$(dS_3,cg_{dS_3})$, or with $(AdS_3,cg_{AdS_3})$; or
$\mathfrak{g}=\mathfrak{so}(3)$, and $(M,g)$ is locally isometric
to the product of $(\mathbb{R},-(dt)^2)$ either with
$(S^3,cg_{S^3})$, or with $(L^3,cg_{L^3})$.
\end{description} Here $c>0$ is a constant, and $S^n$, $L^n$, $dS_n$, $AdS_n$ denote, respectively the sphere, Lobachevskian space, de Sitter space and anti de Sitter space with there
standard metrics. The standard Friedmann-Robertson-Walker
spacetimes are conformally flat and give examples of holonomy
$\mathfrak{so}(1,3)$ \cite{H-L}.

Possible holonomy algebras of conformally flat 4-dimensional
Lorentzian manifolds are classified  also in \cite{H-L}. The first
metric from Theorem \ref{ThC2} in dimension 4 is given in
\cite{St}. In \cite{H-L}, it is stated that it is an open problem
to construct a conformally flat metric with the holonomy algebra
$\mathfrak{sim}(2)$ (which is denoted in \cite{H-L} by $R_{14}$).
An attempt to construct such metric is made in \cite{GT}, where
the following metric is constructed:
\begin{equation}\label{GT01B} g=2dxdt+4ydtdy-4zdtdz+\frac{(dy)^2}{2y^2}+\frac{(dz)^2}{2y^2}+2(x+y^2-z^2)^2(dt)^2.\end{equation}
This metric is conformally flat. Making the transformation
$$x\mapsto x-y^2+z^2,\quad y\mapsto y, \quad z\mapsto z, \quad
t\mapsto t,$$ we obtain
\begin{equation} g=2dxdt+2x^2(dt)^2+\frac{(dy)^2}{2y^2}+\frac{(dz)^2}{2y^2}.\end{equation}
We get that the metric \eqref{GT01B} is decomposible and its
holonomy algebra coincides with
$\mathfrak{so}(1,1)\oplus\mathfrak{so}(2)$, but not with
$\mathfrak{sim}(2)$. Thus in this paper we get metrics with the
holonomy algebra $\mathfrak{sim}(2)$ for the first time (even
more, recall that we find all such metrics).


\end{document}